\documentclass[reqno]{amsart}

\usepackage{amsmath,amsfonts,amssymb,latexsym}
\usepackage[all]{xy}
\usepackage{hyperref}
\usepackage{graphicx, epsfig}
\usepackage{times}
\usepackage{euscript}
\usepackage{ccaption}
\newfixedcaption{\figcaption}{figure}

\def\QQ{\mathbb{Q}}
\def\ZZ{\mathbb{Z}}
\def\RR{\mathbb{R}}

\def\bcp{\mathbb C\mathbb P}
\def\cpb{\overline{\mathbb C \mathbb P ^2}}

\newtheorem{thm}{Theorem}[section]
\newtheorem{prop}[thm]{Proposition}

\newtheorem{cor}[thm]{Corollary}

\newtheorem{rmk}[thm]{Remark}

\numberwithin{equation}{section}

\makeatother

\begin{document}

\title{Smooth structures and Einstein metrics on ${\bcp ^2}\# 5,6,7{\overline{\bcp^2}}$ }

\author{Rare{\c s} R{\u a}sdeaconu}

\address{Rare{\c s} R{\u a}sdeaconu,
         IRMA, 7 Rue Ren{\'e} Descartes, 67184 Strasbourg Cedex,
         FRANCE}

\email{rasdeaco@math.u-strasbg.fr}

\author{Ioana {\c S}uvaina}

\address{Ioana {\c S}uvaina, 
         IHES,
         35 Rue de Chartres,
Bures-sur-Yvette, 91440, FRANCE}
        
\email{ioana@ihes.fr}



\date{\today}

\begin{abstract}
We show that each of the topological 4-manifolds $\bcp^2\# k\overline{\bcp^2},$ for $k=6,7$ admits a smooth structure which has an Einstein metric of scalar curvature $s>0$, a
smooth structure which has an Einstein metric with $s<0$ and
infinitely many non-diffeomorphic smooth structures which do not
admit Einstein metrics. We show that there are infinitely many manifolds homeomorphic  non-diffeomorphic to $\bcp^2\# 5\cpb$ which don't admit an Einstein metric. We also exhibit new examples of manifolds carrying Einstein metrics of both positive and negative scalar curvature.
The main ingredients are recent constructions of exotic symplectic or complex manifolds with small topological numbers.
\end{abstract}

\maketitle

\section{Introduction}
\label{intro}

Recently, new methods were developed \cite{park, abp} to construct symplectic 4-mani-folds with small topology and exotic smooth structures. Moreover, the method proposed by J. Park \cite{park} was later refined \cite{lee-park, ppsh} to produce interesting examples of minimal complex surfaces of general type. In this paper we show how we can use these constructions in regard to the existence or non-existence of Einstein metrics.

 LeBrun and Catanese \cite{leca} showed that there exist homeomorphic, but non-diffeo-morphic manifolds such that one of them admits an Einstein metric of positive sign while the other admits an Einstein metric of negative sign.
The manifold which has positive scalar curvature Einstein metric is $\bcp^2\#8\cpb,$ while the second manifold is a deformation of the Barlow surface.  The key ingredient in their proof was to show that the Barlow surface, a simply connected minimal surface of general type with $K^{2}=1$ can be deformed to a surface with ample canonical bundle. In particular, they showed that there exists a simply connected complex surface of general type with $p_{g}=0, K^{2}=1$ and having ample canonical bundle. The ampleness of the canonical bundle ensures the existence of a K\"ahler-Einstein metric of negative scalar curvature.

Based on ideas from \cite{park}, Lee and Park \cite{lee-park} and more recently Park, Park, and Shin \cite{ppsh} constructed new examples of simply connected, minimal surfaces of general type with 
$p_{g}=0, K^{2}=1,2$ or $3.$ We show that their examples satisfy the ampleness condition:


\begin{thm}\label{amplu}
There exist simply connected complex surfaces of general type, with $p_g=0,~~K^2=2$ or $3,$ and ample canonical bundle.
\end{thm}

In his paper, Park \cite{park}, used his examples and LeBrun's Obstruction Theorem to exhibit a manifold homeomorphic, non-diffeomorphic  to $\bcp^2\#8\cpb$  which does not admit an Einstein metric. We extend his results to $\bcp^2\#k\cpb, k=6,7,$ and we also exhibit an infinite family of homeomorphic manifolds for which no Einstein metric exists.

\begin{thm}
\label{einstein}
Each of the topological 4-manifolds $\bcp^2\# k\overline{\bcp^2},$ for $k=6,7,8$ admits a smooth
structure which has an Einstein metric of scalar curvature $s>0$, a
smooth structure which has an Einstein metric with $s<0$ and
infinitely many non-diffeomorphic smooth structures which do not
admit Einstein metrics.
\end{thm} 

We can also prove non-existence results for an even smaller topological manifold:
\begin{prop}
\label{5-blow-ups}
The canonical smooth structure on $M=\bcp^2\# 5\cpb$ admits an Einstein metric of positive scalar curvature, but $M$ also supports infinitely many exotic smooth structures which don't have an Einstein metric.
\end{prop}

It is a well known fact that in dimensions $2$ and $3$ the sign of the Einstein metric is a topological invariant. This led Besse \cite[p.19]{bes} to consider the conjecture that no smooth compact $n-$manifold can admit Einstein metrics with different scalar curvature signs. As we see in the previous theorems, in dimension $4,$ there are homeomorphic, non-diffeomorphic $4-$manifolds which admit Einstein metrics of opposite sign. For higher dimensions, as far as the authors know, LeBrun and Catanese constructed the only known counterexamples to Besse's conjecture. These examples are in dimensions $4k,$ where $k\geq 2.$ In these dimensions, we are able to provide many new examples with the same property:

\begin{prop}\label{high}
Let $N_1=\bcp^2\# 8\cpb, N_2=\bcp^2\# 7\cpb, N_3=\bcp^2\# 6\cpb.$ Then the manifold $N$ obtained by taking the  $k-$fold products, $k\geq 2,$ of arbitrary  $N_1, N_2,$ or $N_3$, admits two Einstein metrics $g_1,g_2$ such that the signs of the scalar curvature are $s_{g_1}=-1, s_{g_2}=+1.$ Moreover, these metrics are K\"ahler-Einstein with respect to two distinct complex structues $J_1,J_2.$
\end{prop}

As a consequence of the K\"ahlerian property of the metrics and of our method of construction, we have the following in peculiar fact: their corresponding volumes are equal $Vol_{g_1}(N)=Vol_{g_2}(N).$

\begin{rmk}
\label{rmk}
Soon after this paper was finished, Jongil Park and his collaborators found \cite{ppsh2} an analogous example of a  complex surface of general type homeomorphic to $\bcp^2\# 5\cpb.$ Using this construction the results of Theorems \ref{amplu}, \ref{einstein} and Proposition \ref{high} can be extended for $\bcp^2\# 5\cpb.$
 \end{rmk}

The paper is organized as follows: in Section \ref{ampl} we discuss the ampleness of the canonical line bundles of some interesting examples of complex surfaces and the relation to the existence of Einstein metrics, while in Section \ref{exot} we introduce the Seiberg-Witten invariants and use them to get results on non-existence of Einstein metrics. In Section \ref{High} we treat the higher dimensional case.

\section{Ampleness of the (anti)canonical bundle: Existence of Einstein metrics}
\label{ampl}

In general the existence or non-existence of Einstein metrics on a given manifold is hard to prove.
In the case $c_1(M)>0,$ in the unobstructed situations, i.e. for certain complex surfaces whose underlying differential structure is $\bcp^2\#k\cpb, 3\leq k\leq8,$ the existence of K\"ahler-Einstein metrics was proved by Siu \cite{siu}, or Tian and Yau \cite{tian-yau}. A complete solution to the existence was given by Tian \cite{tian}:

\begin{thm}[Tian]\label{tian} 
A compact complex surface $(M^4,J)$ admits
a compatible K\"ahler-Einstein metric with
$s >  0$ if and only if
its anti-canonical line bundle $K_M^{-1}$ is ample and its
Lie algebra of holomorphic vector fields is
reductive. \end{thm}

In the case when $c_1(M)<0,$  a criterion for the existence of a K\"ahler-Einstein metric was independently found by Aubin \cite{aubin} and Yau \cite{yau}:

\begin{thm}[Aubin-Yau]\label{ay}
A compact complex manifold $(M^4,J)$ admits
a compatible K\"ahler-Einstein metric with
$s < 0$ if and only if its
canonical line bundle $K_{M}$ is ample. When such a metric exists, it is 
unique, up to an overall multiplicative constant. 
\end{thm}

In this section, we discuss some examples of smooth structures on  $\bcp^2\#k\cpb, 6\leq k\leq8.$ For each such $k,$ the natural smooth structure as the complex projective plane blown-up at $k$ generic points, can be endowed \cite{tian} by a K\"ahler-Einstein metric of positive scalar curvature. On the other hand, Lee and Park \cite{lee-park} and more recently Park, Park and Shin \cite{ppsh} constructed new exotic smooth structures on $\bcp^2\#k\cpb,$ for any  $6\leq k\leq8.$ Moreover, these admit complex structures yielding interesting examples of minimal surfaces of general type. We are going to prove that the canonical bundle of these surfaces is ample, and so they admit K\"ahler-Einstein metrics of negative scalar curvature.

%

The examples constructed in these papers are
very similar. In this article we are going to treat just the more complicated example \cite{ppsh} of a simply connected  minimal surface of general
type with $p_g=0$ and $K^2=3.$ The other examples
follow similarly. We begin by a brief description of the example of a minimal surface of general
type with $p_g=0$ and $K^2=3,$ referring frequently to \cite{ppsh}.

\begin{figure}[hbtb]
\begin{center}
\hspace*{-1 cm}
\setlength{\unitlength}{1mm}
\includegraphics[width=10cm]{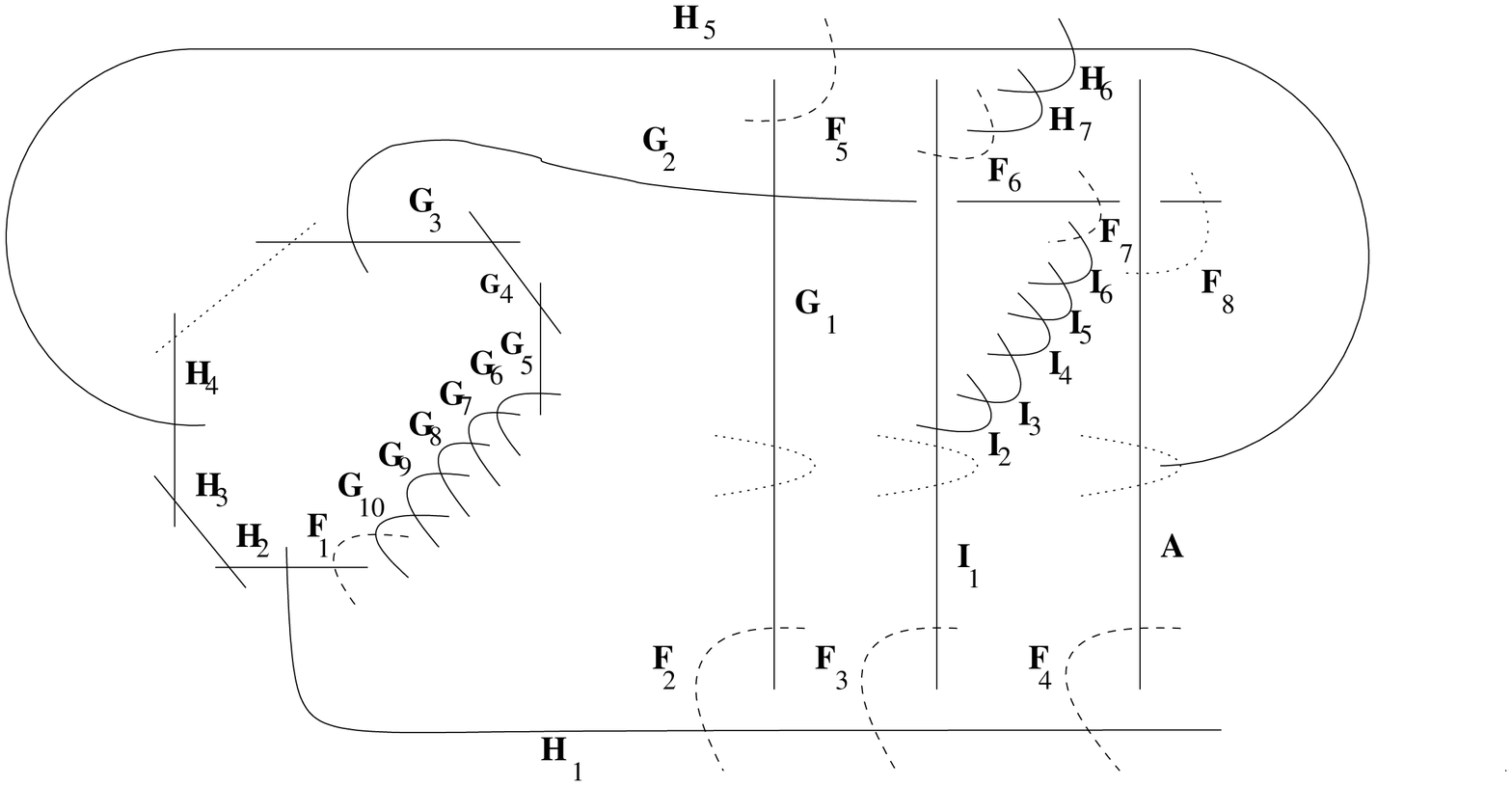}\\
\end{center}
\vspace*{-1 em}
\caption{Manifold $Z=\bcp^{2}\#30\cpb$}
\label{figure}
\end{figure}

We start with an appropriate pencil of plane elliptic curves, and we resolve the base locus to obtain a rational elliptic surface $Y$ with the following set of singular fibers: one $I_8-$fiber, one reducible $I_2-$fiber and two nodal fibers. Moreover, this rational elliptic surface admits 4 sections. Now, $Y$ is blown-up 21 times at well chosen points to obtain a new rational surface $Z.$ These points can be chosen in such a way \cite{ppsh}  that $Z$ contains four particular disjoint chains of rational curves consisting of proper transforms of the singular fibers,  proper transform of some sections and some of the exceptional divisors introduced in the blowing-up process. According to \cite{ppsh}, these chains are $G=\sum_{i=1}^{10} G_i,~H=\sum_{i=1}^7 H_i,~I=\sum_{i=1}^6 I_i$  and a chain of length one denoted by ${A}.$ These chains are represented by continuous lines in Figure \ref{figure}.

They all satisfy Artin's contractibility criterion \cite{km}. By contracting them, we obtain a projective surface $X$ with 4 singularities, each admitting 1-parameter $\QQ-$Gorenstein smoothings. Let $f:Z\to X$ be the contraction map. Now, it is proved \cite{ppsh} that the projective surface $X$  admits a $\QQ-$Gorenstein smoothing, which is a smooth, simply-connected, minimal surface of general type with $K_{X}^{2}=3.$ We denote by $F_{i}, i=1,\cdots 7$ the seven smooth rational curves of self-intersection $-1,$ which are marked by discontinuous lines in Figure \ref{figure} and by $F_8$ the remaining exceptional divisor from the last singular fiber, drawn as a dotted line .  It is easy to see that the Poincar{\'e} duals of the curves $F_{i},i=1,\dots,7,$ together with all of the irreducible components of our 4 chains $G, H, I$ and $A$
form a basis of $H^{2}(Z,\QQ),$ which is torsion-free and 31-dimensional.
%

First, we prove the ampleness of the canonical divisor of the singular surface $X.$ Note that since $X$ has only 1-parameter $\QQ-$Gorenstein smoothing singularities, its canonical divisor is in fact $\QQ-$Cartier.

A common feature of all of the examples appearing in \cite{lee-park, ppsh} is that the pullback of the canonical divisor of the singular variety $X$ to its minimal resolution is effective. For instance, in the example described above, we can write:
\begin{equation}
\label{simplecan}
f^*K_X\equiv_{\QQ} \sum_{i=1}^8 a_iF_i +\sum_{i=1}^{10} b_iG_i + \sum_{i=1}^7 c_iH_i +
\sum_{i=1}^6 d_i I_i + e{A},
\end{equation}
In \cite{ppsh}, the coefficients are explicitely given, but for our purpose it is important to know only that
$a_i, b_i, c_i, d_i$ and $e$ are strictly positive rational numbers. For convenience, let
$Exc(f)=\sum G_i +\sum H_i+\sum I_i+{A}$ denote the exceptional divisor of $f.$

By a direct computation, in \cite{ppsh} it is proved that $K_X^2=3$ and the nefness of $K_X.$
Suppose now that $K_X$ fails to be ample. Since $K_X$ is already nef, according to the Nakai-Moishezon criterion \cite{km}, there exists an irreducible curve $C\subset X$ such that 
$(K_X\cdot C)=0.$ 

The total transform of $C$ in $Z$ is
\begin{eqnarray}
f^*C\equiv_{\QQ} C'+\sum_{i=1}^{10} x_iG_i + \sum_{i=1}^7 y_iH_i + 
\sum_{i=1}^6 z_i I_i + t{A},
\end{eqnarray}
 where $C'$ is the strict transform and $x_i, y_i, z_i, t$ are non-negative 
rational numbers. We should note here that $C'$ is not numerically equivalent to $0.$ To give an immediate proof for this assertion, note that $(C\cdot H)>0$ for any ample line bundle $H$ on 
the projective surface $X.$ But, if $C'\equiv_{\QQ} 0$ we would have 
$$
0=(C'\cdot f^*H)=(f^*C\cdot f^*H)=(C\cdot H),
$$
and this is a contradiction.

Now, a straightforward computation gives:
\begin{align}
(K_X\cdot C)= &(f^*K_X\cdot f^*C)=(f^*K_X\cdot C') \\ \notag 
=&\sum_{i=1}^8 a_i(F_i\cdot C')  +\sum_{i=1}^{10} b_i(G_i\cdot C') + \sum_{i=1}^7 c_i(H_i\cdot C') \\ \notag
&+\sum_{i=1}^6 d_i(I_i\cdot C') + e({A}\cdot C').
\end{align}
Now, obviously the intersection number of $C'$ with any component of $Exc(f)$ is greater or equal to $0,$ with equality if and only if $C'$ does not intersect any of the irreducible components of $Exc(f),$ i.e. $C$ does not pass through the singular points of $X.$ Hence, it follows that 
$$\sum_{i=1}^{10} b_i(G_i\cdot C') + \sum_{i=1}^7 c_i(H_i\cdot C') +\sum_{i=1}^6 d_i(I_i\cdot C') + e({A}\cdot C')\geq 0.
$$ 

Thus $\sum_{i=1}^8 a_i(F_i\cdot C')\leq 0.$ In this case, either there is an $i_0\in \{1,...,8\}$  such that $(C'\cdot F_{i_0})<0,$   or $(C'\cdot F_i)=0$ for all $i=1,...,8,$ and $C'$ does not meet any component of $Exc(f).$ In the first case $C'$ must coincide with $F_{i_0}.$ But the computations in \cite{ppsh,lee-park} show that $(f^*K_X\cdot F_i)>0$ for all $i=1,...,8,$ which is impossible by our assumption. Therefore $C'$ must have vanishing intersection number with all of the $F_i$'s and also with the all of the components of $Exc(f).$ But, as the Poincar\'e duals of the irreducible components of $Exc(f)$ and of the $F_i'$s generate $H^2(Z, \QQ),~C'$ must be numerically trivial on $Z,$ and this is a contradiction. Hence, $X$ has ample canonical bundle.

As a consequence of  Theorem \ref{ay} we get:
\begin{cor}\label{ample}
There exist  smooth complex structures on $\bcp^2\# k\overline{\bcp^2},$ for $k=6,7,8$ which have ample canonical line bundles and hence admitting
K\"ahler-Einstein metrics of negative scalar curvature.
\end{cor} 
\begin{proof}
As we have already mentioned earlier, we only discuss the case of $\bcp^2\# 6\overline{\bcp^2}.$ The other cases follow similarly. 

We know that a  $\QQ-$Gorenstein smoothing of $X$ exists \cite{ppsh} and gives a smooth exotic structure on  $\bcp^2\# 6\overline{\bcp^2}.$ 
As we showed, the singular projective surface $X$ has ample $\QQ-$ Cartier canonical divisor.
But, ampleness is an open property \cite{km}, so the canonical bundle of these smoothings is automatically ample. This will allow us now to apply Aubin-Yau's Theorem \ref{ay} to prove the existence of K\"ahler-Einstein metrics of negative scalar curvature on the minimal surfaces of general type described above.
\end{proof}

\begin{rmk}
\label{barlow}
The examples of complex surfaces with ample canonical bundle we found are new to our knowledge, when $K^{2}=2$ or $3.$ When $K^{2}=1,$ LeBrun and Catanese \cite{leca} were able to find a complex surface with ample canonical bundle homeomorphic to  $\bcp^2\# 8\overline{\bcp^2},$ as a deformation of the Barlow surface. As was pointed out by \cite{lee-park}, it is not known whether the surfaces we found in Corollary \ref{ample} are deformation equivalent to the Barlow surface, or even diffeomorphic to it.
\end{rmk}

\section{Exotic structures: Non-Existence of Einstein metrics}
\label{exot}
The classical obstruction to the existence of an Einstein metric on an oriented, compact, smooth, Riemannian, four-manifold $(M,g)$ is the Hitchin-Thorpe Inequality:
\begin{equation}\label{ht}
(2\chi\pm 3\tau)(M) \geq 0  
\end{equation}
with equality when the Einstein manifold $(M,g)$ is covered by the hyperk\"ahler $K3$ surface or by the flat four-torus, or by the above with the reversed orientations. Here $\chi(M), \tau(M)$ denote the Euler characteristic and the signature of the manifold $M,$ respectively.
Using the Seiberg-Witten theory, LeBrun \cite{lno,lric} was able to find new obstruction to the existence of Einstein metrics. The novelty of his results is that the existence of Einstein metrics depends on the differential structure considered. The techniques developed by LeBrun rely on the existence of a non-trivial solution of the Seiberg-Witten equations.  We briefly introduce the needed notions of the Seiberg-Witten theory and state the main obstruction theorem.

For simplicity, we consider only the case when $H_1(M,\ZZ)$ has no $2-$torsion. Then there is a one-to-one correspondence between the set of $Spin^c$ structures and the set $\{c\in H^2(M,\ZZ)| c_1\equiv w_2(M) \mod 2\}$ of characteristic classes in $H^2(M,\ZZ).$ The Seiberg-Witten invariant is defined \cite{mo} as an integer valued function 
\begin{equation}
SW_M:\{c\in H^2(M,\ZZ)~|~~ c\equiv w_2(M)\mod 2\} \longrightarrow\ZZ  \notag
\end{equation}
We call $c$ a basic class if $SW_M(c)\neq 0.$

The Seiberg-Witten invariant is a diffeomorphism invariant in the case $b_+>1.$ In the case of $b_+=1,$ the invariant depends on the choice of an orientation of $H^2(M,\ZZ)$ and $H^1(M,\RR).$ But, as we are interested in manifolds for which the existence of an Einstein metric is topologically unobstructed,  the inequality (\ref{ht}) must be strictly satisfied.
Hence, in the case when the $SW_M(c)\neq0,$ we have $c^2\geq (2\chi+3\tau)(M) >0.$ But for any arbitrary metric $g,~c^2=c_+^2-c_-^2>0$ where we denote by $c_{\pm}$ the (anti)-self-dual part of the harmonic representative 2-form. Thus $c_+\neq0,$ and in this special situation  (see \cite[Thm 6.9.2]{mo} ) the Seiberg-Witten invariant does not depend on the choice of metric. 

We are now ready to state LeBrun's obstruction, see \cite[Thm. 3.3]{lric}:

\begin{thm}[LeBrun]\label{obstr}
Let $X$ be a compact oriented $4-$manifold with a non-trivial Seiberg-Witten invariant and with $(2\chi+3\tau)(X)>0.$ Then
\begin{equation*}
M=X\#k\cpb 			
\end{equation*}
does not admit an Einstein metric if $k\geq \frac13(2\chi+3\tau)(X).$
\end{thm}

In the last years there was a lot of work done on constructing exotic structures with small topology. In this paper we will use a construction due to Akhmedov, Baykur and Park \cite{abp}. Similar results could also be found in papers of above authors or of Baldridge and Kirk.

Here we are interested in some exotic structures on $M=\bcp^2\#3\cpb.$ In \cite[Section 3.4]{abp}, the authors construct a minimal symplectic manifold $X$ homeomorphic to $M,$ but not diffeomorphic. To give the reader a better understanding of this new manifold, we sketch its construction. We start with two product  manifolds $Y=T^2\times \Sigma_2$ and $T^4$, where  $T^n$ is a torus of dimension $n$ and $\Sigma _2$ is a Riemann surface of genus $2.$ If we consider the diagonal class $[T^2\times point] + [point\times T^2]\in H_2(T^4,\ZZ),$ then a smooth representative in this class is a Riemann surface of genus $2$ with self-intersection $2.$ Blowing up two points on the smooth representative gives $\Sigma'\subset (T^4\#2\cpb),$ a smooth Riemann surface of genus $2$ and $0-$self-intersection. If we endow the manifolds $Y$ and $T^4\#2\cpb$ with the canonical symplectic structures, we can choose the representative $\Sigma'$ to be a symplectic submanifold. We can now take the symplectic sum $X'=Y\#_{\Sigma_2}(T^4\#2\cpb)$ along the two genus $2$ symplectic submanifolds $point\times \Sigma_2 \subset Y$ and $\Sigma' \subset (T^4\#2\cpb).$ $X'$ is a symplectic manifold with topological invariants: Euler characteristic $\chi(X')=6$ and signature $\tau(X')=-2.$ As $X'$ is not simply connected we need to make six Luttinger surgeries along well chosen tori (see \cite{abp}), which yield a new simply connected, minimal, symplectic manifold $X,$ with the same topological invariants as $X'.$
Notice that the definition of the Seiberg-Witten invariant in \cite{abp} is for the Poincar\'e dual of the characteristic class. The results are nevertheless the same.
The manifold $X$ has a unique basic class $\beta=c_1(K_X)\in H_2(X,\ZZ),$ such that its Seiberg-Witten invariant $SW_X(\beta)=1.$ Moreover, $X$ contains a null-homologous torus $\Lambda,$ with a preferred simple loop $\lambda.$ 
In Section 4, \cite{abp}, it is shown that $0-$surgery on $\Lambda$ with respect to $\lambda$ yields a symplectic manifold $X_0,$ which also has an unique basic class $\beta_0=c_1(K_{X_0}).$ Then $1/n-$surgery$,n\geq1$ on $\Lambda$ with respect to $\lambda$  generates a family of manifolds $X_n$ which are homeomorphic to $X_1=X.$ Corresponding to the basic class $\beta$ of $X$ there is a unique basic class $\beta _n$ of $X_n$ for which the Seiberg-Witten invariant is non-zero,  and it can be computed as follows:
$$SW_{X_n}(\beta _n)=SW_X(\beta)+(n-1)SW_{X_0}(\beta _0)=1+(n-1)=n.$$
Hence, the family $X_n, n\geq1$ consists of homeomorphic, pairwise non-diffeomorphic manifolds. 
As $\beta_n$ is the unique basic class on $X_n,$ and the Seiberg-Witten invariant of $X_n, n>1$ is not $\pm1,$ the manifold $X_n$ does not admit a symplectic structure. Of course, its basic class satisfies  $\beta_n^2=\beta^2=6=(2\chi+3\tau)(X)=(2\chi+3\tau)(X_n).$

Let $M_{n,k}=X_n\#k\cpb, k=2,3,4,5.$ Then, the conditions of Theorem \ref{obstr} are satisfied and we have the following:
\begin{prop}\label{non-Einstein}
The manifolds $\bcp^2\#l\cpb$ for $l=5,6,7,8$ support infinitely many non-diffeomorphic exotic smooth structures none of which admits an Einstein metric.
\end{prop}
\begin{proof}
The only thing that needs to be argued is that the manifolds $\{M_{n,k}\}_n$ are homeomorphic but pairwise not diffeomorphic. For this we use again the Seiberg-Witten invariant and the commutativity of our construction. First we blow-up $k=l-3$ points on $X=X_1$ and denote by $E_1,\dots E_{k}$ the exceptional divisors.  This gives a new symplectic manifold, which has basic classes $\beta\pm c_1(E_1)\pm \dots\pm c_1(E_{k}).$ $0-$surgery on $\Lambda$ with respect to $\lambda$ gives the manifold $X_0\#k\cpb$ with corresponding basic classes $\beta_0\pm c_1(E_1)\pm \dots\pm c_1(E_{k}).$ We have a similar relation for the Seiberg-Witten invariant of $M_{n,k}$:
\begin{align}
SW_{M_{n,k}}(\beta _n\pm c_1(E_1)\pm \dots\pm c_1(E_{k}))&=SW_X(\beta\pm c_1(E_1)\pm \dots\pm c_1(E_{k})) \notag \\
&+(n-1)SW_{X_0}(\beta _0\pm c_1(E_1)\pm \dots\pm c_1(E_{k}))\notag \\
&=1+(n-1)=n. \notag
\end {align}
Hence, for any fixed $k\in\{2,3,4,5\},$ the manifolds $M_{n,k}$ are homeomorphic to $\bcp^2\#(k+3)\cpb,$ but pairwise not diffeomorphic. Moreover, if $n>1$ they don't support a symplectic structure.
\end {proof}

\begin{rmk}
We would like to note that since for all of the above smooth structures the Seiberg-Witten invariant is non-trivial the manifolds don't admit Riemannian metrics of positive constant scalar curvature. Moreover, as the $c_1^2(M_{n,k})>0$ there is a bound \cite{lno} on the scalar curvature:
$$\int_{M_{n,k}} s_g^2 d\mu \geq 32 \pi^2c_1^2({M_{n,k}})>0$$
where $d\mu$ is the volume form with respect to $g.$ Hence the manifolds $M_{n,k}$ don't admit non-negative constant scalar curvature metrics, and their Yamabe invariant is negative.
\end{rmk}

The results of Theorem \ref{einstein} are immediate by putting together the results from \ref{ay}, \ref{ample}, \ref{non-Einstein}.

\section{Higher Dimensional Manifolds}\label{High}

Let $S_{1}$ be a complex surface with ample canonical bundle and $c_{1}^{2}=1.$ As an example of such surface, we can either take the one provided by Corollary \ref{ample}, or the example \cite{leca} constructed by LeBrun and Catanese. Let $S_2, S_3$ also be complex surfaces of general type, with ample canonical line bundle and $c_1^2=2,3,$ respectively, as in Corollary \ref{ampl}. We are now ready to give the proof of Proposition \ref{high}.


\begin{proof}
Let $N=N_{i_1}\times\dots\times N_{i_k}$ and let $S=S_{i_1}\times\dots\times S_{i_k},$ where $i_k$ is either $1,2$ or $3.$ The manifolds $S_i$ and $N_i$ are homeomorphic, hence by a theorem of Wall \cite{wall}, they are h-cobordant. These h-cobordisms induce an h-cobordism between $N$ and $S.$ But as $N,S$ are simply connected manifolds of dimension greater than $5$, Smale's h-cobordism theorem tells us that they are diffeomorphic. We know that $N_i$ and $S_i$ admit K\"ahler-Einstein metrics of positive, negative scalar curvature, respectively. We can rescale these metrics  such that the scalar curvatures are $\pm1.$ On $N,S$ we consider the product metrics associated to the corresponding $i_k.$ As we take products of K\"ahler-Einstein with the same scalar curvature, the new metrics are going to be K\"ahler-Einstein. The product complex structures on $N, S$ are of Kodaira dimension $-\infty$ and $4k,$ respectively.
\end{proof}

\begin{rmk}
As the dimension considered increases, the number of manifolds constructed in Proposition \ref{high}  is 
$\frac12 (k+1)(k+2).$\end{rmk}

\section*{Acknowledgements} The authors would like to thank  Claude LeBrun for suggesting this problem and for his remarks and suggetions on an early version of this paper. This paper was written while the first author was a CNRS post-doc at IRMA, Strasbourg, and the second was visiting IHES as a IPDE fellow. We would like to thank these institutions for their hospitality.

\bibliographystyle{alpha}

\end{document}